\documentclass[a4paper]{article}
\usepackage{amssymb}
\usepackage{amsmath}

\title{5-TH ORDER DIFFERENTIAL EQUATIONS RELATED TO CALABI-YAU DIFFERENTIAL
EQUATIONS}

\author{Gert Almkvist}
\begin{document}

\maketitle

\textbf{Introduction.}

The paper [2] started with the 5-th order differential equation \#32 ( where $\theta =x\frac{d}{dx}$ )%
\begin{equation*}
\theta ^{5}-3x(2\theta +1)(3\theta ^{2}+3\theta +1)(15\theta ^{2}+15\theta
+4)-3x^{2}(\theta +1)^{3}(3\theta +2)(3\theta +4)
\end{equation*}%
It came from number theory (approximating $\varsigma (4)$ a l\'{a} Ap\'{e}ry
), see Zudilin [10]. Let 
\begin{equation*}
w_{0}=\sum_{n=0}^{\infty }A_{n}x^{n}
\end{equation*}%
be the analytic solution of \#32. Then 
\begin{equation*}
(n+1)^{5}A_{n+1}-3(n+1)(3n^{2}+3n+1)(15n^{2}+15n+4)A_{n}-3n^{3}(3n-1)(3n-1)A_{n-1}=0
\end{equation*}%
Let $\left\{ A_{n}\right\} $ and $\left\{ B_{n}\right\} $ be two solutions
of the recursion with initial values%
\begin{equation*}
A_{0}=1,A_{1}=12
\end{equation*}%
\begin{equation*}
B_{0}=0,B_{1}=13
\end{equation*}%
Then%
\begin{equation*}
\frac{B_{n}}{A_{n}}\rightarrow \varsigma (4)=\frac{\pi ^{4}}{90}
\end{equation*}%
as $n\rightarrow \infty .$ Krattenthaler and Rivoal found the explicit
formula%
\begin{equation*}
A_{n}=\sum_{i,j}\binom{n}{i}^{2}\binom{n}{j}^{2}\binom{n+i}{n}\binom{n+j}{n}%
\binom{i+j}{n}
\end{equation*}%
(see [5] and[10] for several such formulas). Later Zudilin [10] found a
formula for $A_{n}$ as a simple harmonic sum 
\begin{equation*}
A_{n}=\sum_{k=0}^{n}\binom{n}{k}^{4}\binom{n+k}{n}^{2}\binom{2n-k}{n}%
^{2}\left\{ 1+k(-6H_{k}+6H_{n-k}+2H_{n+k}-2H_{2n-k})\right\} 
\end{equation*}%
where 
\begin{equation*}
H_{k}=\sum_{j=1}^{k}\frac{1}{j}\text{ if }k>0\text{ and }0\text{ otherwise}
\end{equation*}%
The connection to the 4-th order Calabi-Yau differential equations is the
following: Given the 4-th order equation%
\begin{equation*}
y^{(4)}+a_{3}y^{(3)}+a_{2}y^{\prime \prime }+a_{1}y^{\prime }+a_{0}y=0
\end{equation*}%
which is MUM (Maximal Unipotent Monodromy) i.e. it has a Frobenius basis of
solutions%
\begin{equation*}
y_{0}=1+A_{1}x+A_{2}x^{2}+...
\end{equation*}%
\begin{equation*}
y_{1}=y_{0}\log (x)+B_{1}x+...
\end{equation*}%
\begin{equation*}
y_{2}=\frac{1}{2}y_{0}\log (x)^{2}+(B_{1}x+...)\log (x)+C_{1}x+...
\end{equation*}%
\begin{equation*}
y_{3}=\frac{1}{6}y_{0}\log (x)^{3}+\frac{1}{2}(B_{1}x+...)\log
(x)^{2}+(C_{1}x+...)\log (x)+D_{1}x+...
\end{equation*}%
We also assume that the differential equation satisfies the Calabi-Yau
condition (C-Y2)%
\begin{equation*}
a_{1}=\frac{1}{2}a_{2}a_{3}-\frac{1}{8}a_{3}^{3}+a_{2}^{\prime }-\frac{3}{4}%
a_{3}a_{3}^{\prime }-\frac{1}{2}a_{3}^{\prime \prime }
\end{equation*}%
This condition is equivalent to that the Wronskian%
\begin{equation*}
w_{0}=%
\begin{vmatrix}
y_{0} & y_{1} \\ 
y_{0}^{\prime } & y_{1}^{\prime }%
\end{vmatrix}%
\end{equation*}%
satisfies a 5-th order differential equation%
\begin{equation*}
w^{(5)}+b_{4}w^{(4)}+b_{3}w^{(3)}+b_{2}w^{\prime \prime }+b_{1}w^{\prime
}+b_{0}w=0
\end{equation*}%
(otherwise the Wronskian satisfies a 6-th order equation). The condition
C-Y2 can be expressed in $b_{2},b_{3},b_{4}$%
\begin{equation*}
b_{2}=\frac{3}{5}b_{3}b_{4}-\frac{4}{25}b_{4}^{3}+\frac{3}{2}b_{3}^{\prime }-%
\frac{6}{5}b_{4}b_{4}^{\prime }-b_{4}^{\prime \prime }
\end{equation*}%
In [4] there are formulas expressing the $a_{i}^{\prime }$ in the $b_{i},$
e.g.%
\begin{equation*}
a_{3}=\frac{2}{x}+\frac{2}{5}b_{4}
\end{equation*}%
(the other formulas are rather long so we do not state them here). One
verifies that \#32 satisfies C-Y2 and one finds the equation stated in "the
big table" [3].

After this Zudilin and I were very excited and thought that this was the way
to find many new C-Y differential equations. But the result was almost the
opposite. Indeed, 5-th order differential equations satisfying C-Y2 are very
rare (except the ones constructed from known 4-th order C-Y equations).

The invariants of importance for a C-Y equation is the Yukawa coupling $K(q)$
and the instanton numbers $N_{k}$ defined by 
\begin{equation*}
q=\exp (\frac{y_{1}}{y_{0}})
\end{equation*}%
\begin{equation*}
K(q)=(q\frac{d}{dq})^{2}(\frac{y_{2}}{y_{0}})=1+\sum_{k=1}^{\infty }\frac{%
k^{3}N_{k}q^{k}}{1-q^{k}}
\end{equation*}%
Two C-Y differential are equivalent if the have the same instanton numbers.
This is the case after the substitution%
\begin{equation*}
y(x)\rightarrow f(x)y(x+c_{2}x^{2}+c_{3}x^{3}+...)
\end{equation*}

Why are there so few 5-th order C-Y equations? One way to find 4-th order
equations is to use "Zeilberger" in Maple on simple sums of 5-fold products
of binomial coefficients. E.g. 
\begin{equation*}
A_{n}=\sum_{k}\binom{n}{k}^{5}
\end{equation*}%
gives us ( after factorization) the 4-th order equation \#22. If we try the
same trick with%
\begin{equation*}
A_{n}=\sum_{k}\binom{n}{k}^{6}
\end{equation*}%
we get a (nonfactorable) equation of order 6 ( and it is not MUM either).
The same is the case with double sums. Actually \#32,44 and 189 are the only
examples of double sums. The other cases we have found contain harmonic sums.

Finally we ask, is it possible to find 5-th order equations where the
coefficients $A_{n}$ are the constant terms of $S_{n}$ where $S$ is a
Laurent polynomial coming from a reflexive polytope in\textbf{\ R}$^{5}$?
\#188 is such a case (see end of the paper).%
\begin{equation*}
\end{equation*}

\textbf{2.1 Yifan Yang's pullback.}

In [4] we discussed in detail Yifan Yang's pullback of a 5-th order equation
satisfying C-Y2 which often cuts the degree of the ordinary pullback in
half.Here we will only illustrate it with the original example \#32.

Given a 5-th order equation satisfying C-Y2 
\begin{equation*}
w^{(5)}+b_{4}w^{(4)}+b_{3}w^{(3)}+b_{2}w^{\prime \prime }+b_{1}w^{\prime
}+b_{0}w=0
\end{equation*}%
the Yifan Yang pullback is%
\begin{equation*}
y^{(4)}+c_{3}y^{(3)}+c_{2}y^{\prime \prime }+c_{1}y^{\prime }+c_{0}y=0
\end{equation*}%
where%
\begin{equation*}
c_{3}=\frac{8}{5}b_{4}
\end{equation*}%
\begin{equation*}
c_{2}=\frac{1}{2}b_{3}+\frac{7}{5}b_{4}^{\prime }+\frac{19}{25}b_{4}^{2}
\end{equation*}%
\begin{equation*}
c_{1}=-\frac{3}{2}b_{2}+\frac{7}{5}b_{3}^{\prime }+\frac{19}{25}b_{3}b_{4}
\end{equation*}%
\begin{equation*}
c_{0}=-\frac{1}{4}b_{1}+\frac{1}{10}b_{2}^{\prime }+\frac{1}{25}b_{2}b_{4}+%
\frac{9}{40}b_{3}^{\prime \prime }+\frac{1}{16}b_{3}^{2}+\frac{1}{25}%
b_{3}b_{4}^{\prime }+\frac{23}{100}b_{3}^{\prime }b_{4}+\frac{9}{250}%
b_{3}b_{4}^{2}
\end{equation*}

\bigskip

In \#32 we have%
\begin{equation*}
b_{4}=\frac{5(81x^{2}+675x-2)}{x(27x^{2}+270x-1)}
\end{equation*}%
\begin{equation*}
b_{3}=\frac{1752x^{2}+11502x-25}{x^{2}(27x^{2}+270x-1)}
\end{equation*}%
\begin{equation*}
b_{2}=\frac{3(804x^{2}+3753x-5)}{x^{3}(27x^{2}+270x-1)}
\end{equation*}%
\begin{equation*}
b_{1}=\frac{816x^{2}+2130x-1}{x^{4}(27x^{2}+270x-1)}
\end{equation*}%
giving%
\begin{equation*}
c_{3}=\frac{8(81x^{2}+675x-2)}{x(27x^{2}+270x-1)}
\end{equation*}%
\begin{equation*}
c_{2}=\frac{26604x^{4}+4428594x^{3}+17854185x^{2}-105732x+149}{%
2x^{2}(27x^{2}+270x-1)^{2}}
\end{equation*}%
\begin{equation*}
c_{1}=\frac{3(122580x^{4}+1831491x^{3}+6322540x^{2}-32634x+37)}{%
x^{3}(27x^{2}+270x-1)^{2}}
\end{equation*}%
\begin{equation*}
c_{0}=\frac{4769856x^{4}+62980416x^{3}+175580928x^{2}-770052x+625}{%
16x^{4}(27x^{2}+270x-1)^{2}}
\end{equation*}

\bigskip

Converting the derivatives to $\theta =x\frac{d}{dx}$ and making the
substitution $\theta \rightarrow \theta -\frac{5}{2}$ we obtain 
\begin{equation*}
\theta ^{4}-x(540(\theta +\frac{1}{2})^{4}+486(\theta +\frac{1}{2})^{2}+%
\frac{57}{4})
\end{equation*}%
\begin{equation*}
+x^{2}(72846(\theta +1)^{4}+\frac{6915}{2}(\theta +1)^{2}+\frac{3}{4})
\end{equation*}%
\begin{equation*}
+x^{3}(14580(\theta +\frac{3}{2})^{4}+12717(\theta +\frac{3}{2})^{2}+324)
\end{equation*}%
\begin{equation*}
+\frac{9}{16}x^{4}(6\theta +11)^{2}(6\theta +13)^{2}
\end{equation*}%
Note that the ordinary pullback in the big table has degree 8. We also
observe that due to the symmetry the differential equation is determined by
only 12 parameters (in the big table there are about 38 parameters!).%
\begin{equation*}
\end{equation*}

\textbf{2.2. 5-th order equations of degree 1}

These are the hypergeometric differential equations of the form%
\begin{equation*}
\theta ^{5}-xP(\theta )
\end{equation*}

where $P(\theta )$ is a polynomial of deqree $5.$This is treated in detail
in section 3.2 in [4]. Here the Yifan Yang pullback is really useful,
resulting in 14 essentially new 4-th order C-Y equations of degree 2.
\bigskip
\textbf{Example. \#}$\widetilde{\mathbf{11}}$

We start with 
\begin{equation*}
\theta ^{5}-4\cdot 12^{3}x(\theta +\frac{1}{2})(\theta +\frac{1}{4})(\theta +%
\frac{3}{4})(\theta +\frac{1}{3})(\theta +\frac{2}{3})
\end{equation*}%
and end with the pullback%
\begin{equation*}
\theta ^{4}-12x(1152\theta ^{4}+2304\theta ^{3}+2710\theta ^{2}+1558\theta
+341)
\end{equation*}%
\begin{equation*}
+144x^{2}(24\theta +19)(24\theta +23)(24\theta +25)(24\theta +29)
\end{equation*}

In [4] we note that the 5-th order equations giving $\widetilde{3},%
\widetilde{4},\widetilde{5},\widetilde{6},\widetilde{8},\widetilde{10},%
\widetilde{11},\widetilde{12},\widetilde{13}$ are not really needed since
these equations are equivalent to certain Hadamard products of two second
order equations. E.g. $\widetilde{11}$ is equivalent to $h\ast i$ which is
of degree 8. But of course we prefer the simple $\widetilde{11}$ before the
clumsy $h\ast i.$%
\begin{equation*}
\end{equation*}

\textbf{2.3. 5-th order equations of degree 2.}

\textbf{(i) }\texttt{\ }Here we have all the Hadamard products $A\ast \alpha
,...,D\ast \varkappa $ (see [2] and [4] for definitions). There are $4\cdot
10-1=39$ of them since the pullback of $B\ast \vartheta $ has trivial $K(q).$
We list the third order equations%
\begin{equation}
\theta ^{3}-x(2\theta +1)(a\theta ^{2}+a\theta +b)+cx^{2}(\theta +1)^{3}
\end{equation}%
\begin{equation*}
\begin{tabular}{|l|l|l|l|}
\hline
Name & a & b & c \\ \hline
$\alpha $ & 10 & 4 & 64 \\ \hline
$\beta $ & 16 & 8 & 256 \\ \hline
$\gamma $ & 17 & 5 & 1 \\ \hline
$\delta $ & 7 & 3 & 81 \\ \hline
$\epsilon $ & 12 & 4 & 16 \\ \hline
$\zeta $ & 9 & 3 & -27 \\ \hline
$\eta $ & 11 & 5 & 125 \\ \hline
$\vartheta $ & 64 & 40 & 4096 \\ \hline
$\iota $ & 27 & 15 & 729 \\ \hline
$\kappa $ & 432 & 312 & 186624 \\ \hline
\end{tabular}%
\end{equation*}

\textbf{Case A.}

The Yifan Yang pullback of%
\begin{equation*}
\theta ^{5}-4x(2\theta +1)^{3}(a\theta ^{2}+a\theta +b)+16cx^{2}(\theta
+1)(2\theta +1)^{2}(2\theta +3)^{2}
\end{equation*}%
is%
\begin{equation*}
\theta ^{4}-x\left\{ 16a(\theta +\frac{1}{2})^{4}+(13a+4b)(\theta +\frac{1}{2%
})^{2}+\frac{a}{4}+\frac{b}{2}\right\}
\end{equation*}%
\begin{equation*}
+x^{2}\left\{ (64a^{2}+32c)(\theta +1)^{4}+(32ab-8a^{2}+108c)(\theta +1)^{2}+%
\frac{a^{2}}{4}+4b^{2}-2ab+\frac{39c}{4}\right\}
\end{equation*}%
\begin{equation*}
-16cx^{3}(\theta +\frac{3}{2})^{2}\left\{ 16a(\theta +\frac{3}{2}%
)^{2}+11a+4b\right\}
\end{equation*}%
\begin{equation*}
+64c^{2}x^{4}(\theta +2)^{2}(2\theta +3)(2\theta +5)
\end{equation*}

\textbf{Case B.}

The Yifan pullback of 
\begin{equation*}
\theta ^{5}-3x(2\theta +1)(3\theta +1)(3\theta +2)(a\theta ^{2}+a\theta
+b)+9cx^{2}(\theta +1)(3\theta +1)(3\theta +2)(3\theta +4)(3\theta +5)
\end{equation*}%
is%
\begin{equation*}
\theta ^{4}-x\left\{ 36a(\theta +\frac{1}{2})^{4}+(29a+9b)(\theta +\frac{1}{2%
})^{2}+\frac{a}{2}+\frac{5b}{4}\right\}
\end{equation*}%
\begin{equation*}
+x^{2}\left\{ (324a^{2}+162c)(\theta +1)^{4}+(162ab-45a^{2}+\frac{217c}{4}%
)(\theta +1)^{2}+\frac{81b^{2}}{4}-9ab+a^{2}+\frac{97c}{2}\right\}
\end{equation*}%
\begin{equation*}
-3^{5}cx^{3}(\theta +\frac{3}{2})^{2}\left\{ 12a(\theta +\frac{3}{2}%
)^{2}+8a+3b\right\}
\end{equation*}%
\begin{equation*}
+2^{4}3^{6}c^{2}x^{4}(\theta +\frac{3}{2})(\theta +\frac{5}{2})(3\theta +%
\frac{11}{2})(3\theta +\frac{13}{2})
\end{equation*}

\textbf{Case C.}.

The Yifan Yang pullback of%
\begin{equation*}
\theta ^{5}-4x(2\theta +1)(4\theta +1)(4\theta +3)(a\theta ^{2}+a\theta
+b)+16cx^{2}(\theta +1)(4\theta +1)(4\theta +3)(4\theta +5)(4\theta +7)
\end{equation*}

is%
\begin{equation*}
\theta ^{4}-x\left\{ 64a(\theta +\frac{1}{2})^{4}+(51a+16b)(\theta +\frac{1}{%
2})^{2}+\frac{3a}{4}+\frac{5b}{4}\right\}
\end{equation*}%
\begin{equation*}
+x^{2}\left\{ (1024a^{2}+512c)(\theta +1)^{4}+(512ab-160a^{2}+1712c)(\theta
+1)^{2}+64b^{2}-24ab+\frac{9a^{2}}{4}+\frac{599c}{4}\right\}
\end{equation*}%
\begin{equation*}
-2^{8}cx^{3}(\theta +\frac{3}{2})^{2}\left\{ 64a(\theta +\frac{3}{2}%
)^{2}+41a+16b\right\}
\end{equation*}%
\begin{equation*}
+2^{12}c^{2}x^{4}(\theta +\frac{3}{2})(\theta +\frac{5}{2})(4\theta
+7)(4\theta +9)
\end{equation*}

\textbf{Case D}

The Yifan Yang pullback of%
\begin{equation*}
\theta ^{5}-12x(2\theta +1)(6\theta +1)(6\theta +5)(a\theta ^{2}+a\theta
+b)+144cx^{2}(\theta +1)(6\theta +1)(6\theta +5)(6\theta +7)(6\theta +11)
\end{equation*}

is%
\begin{equation*}
\theta ^{4}-x\left\{ 144a(\theta +\frac{1}{2})^{4}+(113a+36b)(\theta +\frac{1%
}{2})^{2}+\frac{5a}{4}+\frac{13b}{2}\right\}
\end{equation*}%
\begin{equation*}
+x^{2}\left\{ (5184a^{2}+2592c)(\theta
+1)^{4}+(2592ab-936a^{2}+8604c)(\theta +1)^{2}+324b^{2}-90ab+\frac{25a^{2}}{4%
}+\frac{2927c}{4}\right\}
\end{equation*}%
\begin{equation*}
-2^{4}3^{5}cx^{3}(\theta +\frac{3}{2})^{2}\left\{ 48a(\theta +\frac{3}{2}%
)^{2}+29a+12b\right\}
\end{equation*}%
\begin{equation*}
+2^{6}3^{6}c^{2}x^{4}(2\theta +3)(2\theta +5)(3\theta +5)(3\theta +7)
\end{equation*}%
\begin{equation*}
\end{equation*}

\textbf{(ii) \ }Then we have the "sporadic" equations of degree 2. Firstly
we have \#32 which we already treated.

\begin{equation*}
\end{equation*}
\newpage
\#\textbf{60.}%
\begin{equation*}
A_{n}^{\prime }=\sum_{k}\binom{n}{k}^{2}\binom{n+k}{n}\binom{2n-k}{n}\binom{%
2k}{k}\binom{2n-2k}{n-k}
\end{equation*}%
\begin{equation*}
\theta ^{5}-2x(2\theta +1)(31\theta ^{4}+62\theta ^{3}+54\theta
^{2}+23\theta +4)+12x^{2}(\theta +1)(3\theta +2)(3\theta +4)(4\theta
+3)(4\theta +5)
\end{equation*}%
with Yifan Yang pullback%
\begin{equation*}
\theta ^{4}-x\left\{ 248(\theta +\frac{1}{2})^{4}+232(\theta +\frac{1}{2}%
)^{2}+\frac{15}{2}\right\}
\end{equation*}%
\begin{equation*}
+x^{2}\left\{ 18832(\theta +1)^{4}+13806(\theta +1)^{2}+1251\right\}
\end{equation*}%
\begin{equation*}
-3x^{3}\left\{ 142848(\theta +\frac{3}{2})^{4}+127432(\theta +\frac{3}{2}%
)^{2}+3046\right\}
\end{equation*}%
\begin{equation*}
+9x^{4}(24\theta +41)(24\theta +47)(24\theta +49)(24\theta +55)
\end{equation*}%

\textbf{\#189.}%
\begin{equation*}
A_{n}^{\prime }=\binom{2n}{n}\sum_{i,j}\binom{n}{i}^{2}\binom{n}{j}^{2}%
\binom{i+j}{n}^{2}
\end{equation*}%
\begin{equation*}
\theta ^{5}-2x(2\theta +1)(65\theta ^{4}+130\theta ^{3}+105\theta
^{2}+40\theta +6)+16x^{2}(\theta +1)(2\theta +1)(2\theta +3)(4\theta
+3)(4\theta +5)
\end{equation*}%
with Yifan Yang pullback%
\begin{equation*}
\theta ^{4}-2x\left\{ 260(\theta +\frac{1}{2})^{4}+235(\theta +\frac{1}{2}%
)^{2}+7\right\}
\end{equation*}%
\begin{equation*}
+4x^{2}\left\{ 17412(\theta +1)^{4}+2727(\theta +1)^{2}+181\right\}
\end{equation*}%
\begin{equation*}
-16x^{3}\left\{ 33280(\theta +\frac{3}{2})^{4}+27480(\theta +\frac{3}{2}%
)^{2}+431\right\}
\end{equation*}%
\begin{equation*}
+256x^{4}(8\theta +13)(8\theta +15)(8\theta +17)(8\theta +19)
\end{equation*}%
\begin{equation*}
\end{equation*}

\textbf{\#244.}%
\begin{equation*}
A_{n}^{\prime }=\sum_{k}\binom{n}{k}^{6}\binom{2k}{k}\binom{2n-2k}{n-k}%
\left\{ 1+k(-8H_{k}+8H_{n-k}+2H_{2k}-2H_{2n-2k})\right\}
\end{equation*}%
\begin{equation*}
\theta ^{5}+2x(2\theta +1)(26\theta ^{4}+52\theta ^{3}+44\theta
^{2}+18\theta +3)-12x^{2}(\theta +1)^{3}(6\theta +5)(6\theta +7)
\end{equation*}%
with Yifan Yang pullback%
\begin{equation*}
\theta ^{4}+2x\left\{ 104(\theta +\frac{1}{2})^{4}+96(\theta +\frac{1}{2}%
)^{2}+3\right\}
\end{equation*}%
\begin{equation*}
+x^{2}\left\{ 9952(\theta +1)^{4}-1978(\theta +1)^{2}-309\right\}
\end{equation*}%
\begin{equation*}
-3x^{3}\left\{ 29952(\theta +\frac{3}{2})^{4}+27440(\theta +\frac{3}{2}%
)^{2}+828\right\}
\end{equation*}%
\begin{equation*}
+9x^{4}(12\theta +23)^{2}(12\theta +25)^{2}
\end{equation*}%
\begin{equation*}
\end{equation*}

\textbf{\#245.}%
\begin{equation*}
A_{n}^{\prime }=3\binom{2n}{n}^{2}\sum_{k=0}^{[n/3]}(-1)^{k}\frac{n-2k}{2n-3k%
}\binom{n}{k}^{4}\binom{3n-3k}{2n}\binom{2n}{3k}^{-1}
\end{equation*}%
\begin{equation*}
\theta ^{5}-6x(2\theta +1)(18\theta ^{4}+36\theta ^{3}+34\theta
^{2}+16\theta +3)+2916x^{2}(\theta +1)^{3}(2\theta +1)(2\theta +3
\end{equation*}%
with Yifan Yang pullback%
\begin{equation*}
\theta ^{4}-3x\left\{ 144(\theta +\frac{1}{2})^{4}+140(\theta +\frac{1}{2}%
)^{2}+5\right\}
\end{equation*}%
\begin{equation*}
+9x^{2}\left\{ 7776(\theta +1)^{4}+9918(\theta +1)^{2}+931\right\}
\end{equation*}%
\begin{equation*}
-2^{2}3^{7}x^{3}\left\{ 576(\theta +\frac{3}{2})^{4}+524(\theta +\frac{3}{2}%
)^{2}+12\right\}
\end{equation*}%
\begin{equation*}
+3^{12}x^{4}(4\theta +7)^{2}(4\theta +9)^{2}
\end{equation*}%
Observe that the ordinary pullback also has degree 4.%
\begin{equation*}
\end{equation*}

\#\textbf{253.}%
\begin{equation*}
A_{n}^{\prime }=\sum_{k}\binom{n}{k}^{2}\binom{n+k}{n}\binom{2n-k}{n}\binom{%
2k}{k}^{2}\binom{2n-2k}{n-k}^{2}
\end{equation*}%
\begin{equation*}
\times \left\{
1+k(-7H_{k}+7H_{n-k}+H_{n+k}-H_{2n-k}+4H_{2k}-2H_{2n-2k})\right\} 
\end{equation*}%
\begin{equation*}
\theta ^{5}-4x(2\theta +1)(16\theta ^{4}+32\theta ^{3}+31\theta
^{2}+15\theta +3)+16x^{2}(\theta +1)(4\theta +3)^{2}(4\theta +5)^{2}
\end{equation*}%
with ordinary pullback (the Yifan Yang pullback has degree 5 and lacks
symmetry)

\bigskip

\begin{equation*}
\theta ^{4}-2^{2}x(64\theta ^{4}+32\theta ^{3}+15\theta ^{2}-\theta -2)
\end{equation*}%
\begin{equation*}
+2^{8}x^{2}(96\theta ^{4}+96\theta ^{3}+53\theta ^{2}+13\theta +8)
\end{equation*}%
\begin{equation*}
-2^{12}x^{3}(256\theta ^{4}+384\theta ^{3}+244\theta ^{2}+84\theta +7)
\end{equation*}%
\begin{equation*}
+2^{18}x^{4}(2\theta +1)^{2}(4\theta +1)(4\theta +3)
\end{equation*}%
\newpage

\#\textbf{255. \ \ }%
\begin{equation*}
A_{n}^{\prime }=\sum_{k}\binom{n}{k}^{2}\binom{n+k}{n}\binom{2n-k}{n}\binom{%
2k}{k}\binom{2n-2k}{n-k}\binom{4k}{2k}\binom{4n-4k}{2n-2k}
\end{equation*}%
\begin{equation*}
\times \left\{
1+k(-5H_{k}+5H_{n-k}+H_{n+k}-H_{2n-k}-2H_{2k}+2H_{2n-2k})\right\}
\end{equation*}%
\begin{equation*}
\theta ^{5}+4x(2\theta +1)(64\theta ^{4}+128\theta ^{3}+141\theta
^{2}+77\theta +17)+16x^{2}(\theta +1)(8\theta +5)(8\theta +7)(8\theta
+9)(8\theta +11)
\end{equation*}%
with Yifan Yang pullback

\bigskip

\begin{equation*}
\theta ^{4}+x\left\{ 1024(\theta +\frac{1}{2})^{4}+1076(\theta +\frac{1}{2}%
)^{2}+43\right\}
\end{equation*}%
\begin{equation*}
+2^{7}x^{2}\left\{ 3072(\theta +1)^{4}+4264(\theta +1)^{2}+413\right\}
\end{equation*}%
\begin{equation*}
+2^{12}x^{3}\left\{ 16384(\theta +\frac{3}{2})^{4}+16576(\theta +\frac{3}{2}%
)^{2}+491\right\}
\end{equation*}%
\begin{equation*}
+2^{22}x^{4}(4\theta +7)(4\theta +9)(8\theta +15)(8\theta +17)
\end{equation*}%
\begin{equation*}
\end{equation*}

\#\textbf{281.}%
\begin{equation*}
A_{n}^{\prime }=3\sum_{k=[2(n+1)/3]}^{n}(-1)^{k}\frac{n-2k}{n-3k}\binom{n+k}{%
n}\binom{2n-k}{n}\binom{n}{3n-3k}\binom{3k}{n}^{-1}\frac{(3k)!}{k!^{3}}\frac{%
(3n-3k)!}{(n-k)!^{3}}
\end{equation*}%
\begin{equation*}
\theta ^{5}+x(2\theta +1)(41\theta ^{4}+82\theta ^{3}+74\theta ^{2}+33\theta
+6)+5x^{2}(\theta +1)(5\theta +3)(5\theta +4)(5\theta +6)(5\theta +7
\end{equation*}%
with Yifan Yang pullback%
\begin{equation*}
\theta ^{4}+x\left\{ 164(\theta +\frac{1}{2})^{4}+156(\theta +\frac{1}{2}%
)^{2}+\frac{21}{4}\right\}
\end{equation*}%
\begin{equation*}
+x^{2}\left\{ 12974(\theta +1)^{4}+\frac{45175}{2}(\theta +1)^{2}+\frac{8835%
}{4}\right\}
\end{equation*}%
\begin{equation*}
+5x^{3}\left\{ 102500(\theta +\frac{3}{2})^{4}+92375(\theta +\frac{3}{2}%
)^{2}+2168\right\}
\end{equation*}%
\begin{equation*}
+5^{6}x^{4}(5\theta +\frac{17}{2})(5\theta +\frac{19}{2})(5\theta +\frac{21}{%
2})(5\theta +\frac{23}{2})
\end{equation*}%
\begin{equation*}
\end{equation*}

\textbf{2.4 Fifth order equations of degree 3.}

We have only two equations of degree 3.

\#\textbf{130. (}Verrill)%
\begin{equation*}
A_{n}^{\prime }=\sum_{i+j+k+l+m+s=n}\left( \frac{n!}{i!j!k!l!m!s!}\right)
^{2}
\end{equation*}%
\begin{equation*}
\theta ^{5}-2x(2\theta +1)(14\theta ^{4}+28\theta ^{3}+28\theta
^{2}+14\theta +3)
\end{equation*}%
\begin{equation*}
+4x^{2}(\theta +1)^{3}(196\theta ^{2}+392\theta +255)-1152x^{3}(\theta
+1)^{2}(\theta +2)^{2}(2\theta +3)
\end{equation*}%
with Yifan Yang pullback%
\begin{equation*}
\theta ^{4}-4x\left\{ 28(\theta +\frac{1}{2})^{4}+28(\theta +\frac{1}{2}%
)^{2}+1\right\}
\end{equation*}%
\begin{equation*}
+3x^{2}\left\{ 1568(\theta +1)^{4}+2130(\theta +1)^{2}+225\right\}
\end{equation*}%
\begin{equation*}
-4x^{3}\left\{ 23104(\theta +\frac{3}{2})^{4}+32532(\theta +\frac{3}{2}%
)^{2}+3213\right\}
\end{equation*}%
\begin{equation*}
+x^{4}\left\{ 872704(\theta +2)^{4}+995680(\theta +2)^{2}+93337\right\}
\end{equation*}%
\begin{equation*}
-2^{9}3^{2}x^{5}(\theta +\frac{5}{2})^{2}\left\{ 784(\theta +\frac{5}{2}%
)^{2}+647\right\}
\end{equation*}%
\begin{equation*}
+2^{14}3^{4}x^{6}(2\theta +1)^{2}(\theta +\frac{5}{2})(\theta +\frac{7}{2})
\end{equation*}%
\begin{equation*}
\end{equation*}

\#\textbf{188.}%
\begin{equation*}
A_{n}^{\prime }=\binom{2n}{n}\sum_{i+j+k+l+m=n}\left( \frac{n!}{i!j!k!l!m!}%
\right) ^{2}
\end{equation*}%
\begin{equation*}
\theta ^{5}-2x(2\theta +1)(35\theta ^{4}+70\theta ^{3}+63\theta
^{2}+28\theta +5)
\end{equation*}%
\begin{equation*}
+4x^{2}(\theta +1)(2\theta +1)(2\theta +3)(259\theta ^{2}+518\theta
+285)-1800x^{3}(\theta +1)(\theta +2)(2\theta +1)(2\theta +3)(2\theta +5)
\end{equation*}

with Yifan Yang pullback%
\begin{equation*}
\theta ^{4}-x\left\{ 280(\theta +\frac{1}{2})^{4}+266(\theta +\frac{1}{2}%
)^{2}+9\right\}
\end{equation*}%
\begin{equation*}
+x^{2}\left\{ 27888(\theta +1)^{4}+31638(\theta +1)^{2}+3063\right\}
\end{equation*}%
\begin{equation*}
-x^{3}\left\{ 1189120(\theta +\frac{3}{2})^{4}+1276704(\theta +\frac{3}{2}%
)^{2}+75906\right\}
\end{equation*}%
\begin{equation*}
+x^{4}\left\{ 21204736(\theta +2)^{4}+11869760(\theta +2)^{2}+941769\right\}
\end{equation*}%
\begin{equation*}
-2^{3}15^{2}x^{5}(4\theta +9)(4\theta +11)\left\{ 4144(\theta +\frac{5}{2}%
)^{2}+2435\right\}
\end{equation*}%
\begin{equation*}
+30^{4}x^{6}(4\theta +9)(4\theta +11)(4\theta +13)(4\theta +15)
\end{equation*}

Let 
\begin{equation*}
S=x+y+z+t+u+\frac{1}{x}+\frac{1}{y}+\frac{1}{z}+\frac{1}{t}+\frac{1}{u}
\end{equation*}%
Then 
\begin{equation*}
A_{n}^{\prime }=\text{constant term of }S^{2n}
\end{equation*}%
\begin{equation*}
\end{equation*}

$\bigskip $\textbf{2.5. Degree higher than 3.}

In the big table there are many 5-th order equations of degree higher than
3. But they are all equivalent to equations of degree at most 3 or where we
do not need any 5-th order equation at all. The worst example is \#129 of
degree 6 and the ordinary pullback has degree 24. It is equivalent to $%
\widetilde{8}$ of degree 2. These equations should be removed from the
table. So in fact we do not know of any 5-th order equation of degree larger
than 3 which is needed.

\begin{equation*}
\end{equation*}

\textbf{Acknowledgements.}

I thank the following persons who have made it possible to write this paper:
Christian van Enckevort, Christian Krattenthaler, Duco van Straten, Helena
Verrill, Yifan Yang and Wadim Zudilin.

\begin{equation*}
\end{equation*}

\textbf{References.}

1. G.Almkvist, Str\"{a}ngar i m\aa nsken I, Normat 51 (2003), 22-33, II,
Normat 51 (2003), 63-79.

2. G.Almkvist, W.Zudilin, Differential equations, mirror maps and zeta
values. In Mirror symmetry,

Proceedings of BIRS workshop on Calabi-Yau Varities and Mirror Symmetry,
December 2003.

3. G.Almkvist, C.van Enckevort, D.van Straten, W.Zudilin, Tables of
Calabi-Yau equations, math.AG/0507430

4. G.Almkvist, Calabi-Yau differential equations of degree 2 and 3 and Yifan
Yang's pullback,

math.AG/0612215

5. G.Almkvist, Some binomial identities related to Calabi-Yau differential
equations, math.CO/0703255

6. G.Almkvist, C.Krattenthaler, Some harmonic sums related to Calabi-Yau
differential equations, in

preparation.

7. C.Krattenthaler, T.Rivoal, Hyperg\'{e}om\'{e}trie et fonction z\^{e}ta de
Riemann, math.NT/0311114

8. Yifan Yang, personal communication.

9. D.Zagier, Integral solutions of Ap\'{e}ry-like recurrence equations,
manuscript 2003.

10. W.Zudilin, Binomial sums related to rational approximations of $%
\varsigma (4),$Mat. Zametki 75 (2004), 637-640,

English transl. Math. Notes 75 (2004), 594-597, math.Ca/031196.%
\begin{equation*}
\end{equation*}

Math. Dept., Univ. of Lund

Box 118

22100 Lund, Sweden

gert@maths.lth.se

\end{document}